\documentclass[11pt]{article}
\usepackage{a4}
\usepackage{graphicx}
\usepackage{parskip}

\usepackage{amsfonts}
\usepackage{natbib}
\usepackage{algorithm}
\usepackage{algorithmic}

\usepackage{amsmath, amssymb}

\usepackage{latexsym}

\newtheorem{theorem}{Theorem}[section]

\newtheorem{lemma}{Lemma}[section]
\newtheorem{corollary}{Corollary}[section]

\newtheorem{definition}{Definition}[section]

\newtheorem{condition}{Condition}[section]

\newcommand{\R}{\mathbb{R}}

\newcommand{\PP} {{  \rm I\hskip-0.22em P}}

\newcommand{\EE} {{\rm I\hskip-0.48em E}}

\long\def\comment#1{}

\newcommand{\defn}{\ensuremath{: \, =}}

\newcommand{\pen}{\ensuremath{{\rm{pen}}}}

\newcommand{\GSET}{\ensuremath{\mathcal{G}}}
\newcommand{\dfun}{\ensuremath{\mathrm{d}}}
\newcommand{\calF}{\ensuremath{\mathcal{F}}}
\newcommand{\calG}{\ensuremath{\mathcal{G}}}
\newcommand{\PhiConj}{\ensuremath{\Phi^*}}
\newcommand{\Gconj}{\ensuremath{G^*}}

\newcommand{\QUADE}[1]{\ensuremath{ {\mathcal{Q}}(#1; {\bf E})}}
  \newcommand{\QUADF}[1]{\ensuremath{ {\mathcal{Q}}(#1; {\bf F})}}

\newenvironment{carlist}
 {\begin{list}{$\bullet$}
 {\setlength{\topsep}{0in} \setlength{\partopsep}{0in}
  \setlength{\parsep}{0in} \setlength{\itemsep}{\parskip}
  \setlength{\leftmargin}{0.07in} \setlength{\rightmargin}{0.08in}
  \setlength{\listparindent}{0in} \setlength{\labelwidth}{0.08in}
  \setlength{\labelsep}{0.1in} \setlength{\itemindent}{0in}}}
 {\end{list}}

\newcommand{\bcar}{\begin{carlist}}
\newcommand{\ecar}{\end{carlist}}

\usepackage{natbib}
\usepackage{color}

\begin{document}

\centerline{\bf {\Large{On concentration for (regularized) empirical risk
  minimization}}}

\vskip .1in

\centerline {Sara van de Geer}

\centerline{Seminar for Statistics, ETH Z\"urich}

\vskip .1in
\centerline{and}

\vskip .1in

\centerline {Martin J. Wainwright}

\centerline{Department of Statistics and Department of EECS}

\centerline{University of California, Berkeley}

\vspace*{.1in}

{\bf Abstract.} Rates of convergence for empirical risk minimizers
have been well studied in the literature. In this paper, we aim to
provide a complementary set of results, in particular by showing that
after normalization, the risk of the empirical minimizer concentrates
on a single point. Such results have been established
by~\cite{chatterjee2014new} for constrained estimators in the normal
sequence model. We first generalize and sharpen this result to
regularized least squares with convex penalties, making use of a
``direct'' argument based on Borell's theorem.  We then study
generalizations to other loss functions, including the negative
log-likelihood for exponential families combined with a strictly
convex regularization penalty. The results in this general setting are
based on more ``indirect'' arguments as well as on concentration
inequalities for maxima of empirical processes.


\section{Introduction}
\label{introduction.section}

Empirical risk minimization (ERM) is an important methodology in
statistics and machine learning, widely used for estimating
high-dimensional and/or nonparametric parameters of interest. The idea
is to express the parameter as minimizer of an expected loss, the
so-called population or theoretical risk.  Given that the distribution
of data is not known or difficult to assess, one replaces the
theoretical expectation by an empirical counterpart defined by
samples. The technique of ERM is known under various names, including
$M$-estimation and minimum contrast estimation.

By the law of large numbers, empirical averages of various types of
random variables, with the i.i.d.\ setting being the canonical case,
are close to their expectation. This elementary fact is the motivation
for ERM and the starting point for studying its theoretical
properties.  There is much literature developing the theory for a
broad spectrum of estimation problems. The more recent literature
takes a non-asymptotic point of view, in which context concentration
inequalities play a major role.  Concentration inequalities
describe the amount of concentration of certain (complex) quantities
around their mean. We refer to \cite{Talagrand:95} as a key paper in
the area, and to the important monographs
\cite{ledoux2001concentration} and more recently
\cite{boucheron2013concentration}.  The key point is that the
deviation from the mean is generally of much smaller order than the
mean itself. Moreover, at least in a certain sense, the deviation does
not depend on the complexity of the original object. In statistics, the
usefulness of concentration inequalities has been excellently outlined
and studied in \cite{massart2000some}.  We also refer the reader to
\cite{koltchinskii2011oracle} for an in-depth treatment in the context
of high-dimensional problems.

Some statistical papers address concentration for the parameter of
interest itself; for instance, see \cite{boucheron2011high} and
\cite{saumard2012optimal}.  The present paper is along the lines of
\cite{chatterjee2014new}. The latter examines the concentration
properties of constrained estimators for the normal sequence model, or
alternatively phrased in the regression setting, for the least-squares
problem with fixed design and Gaussian errors.  The author shows that
the statistical error of the least squares estimator satisfies a
concentration inequality where the amount of concentration still
depends on the complexity of the problem, but is in the nonparametric
case of smaller order than the statistical error itself.  In
\cite{muro15}, the situation is studied where a regularization penalty
based on a squared norm is added to the least squares loss
function. In Section \ref{direct.section}, we provide a ``direct''
argument for concentration of the regularized least squares in the
normal sequence setting.  Our argument here is elementary, using
standard facts from convex analysis~\citep{rockafellar1970convex}, and
concentration for Lipschitz functions of Gaussian
vectors~\citep{borell1975brunn}.  Our next goal is to extend such
results to more general problems.  The main obstacle is that the
direct concentration for Lipschitz functions holds only for the
Gaussian case.  Accordingly, we make use of more general one-sided
concentration results for maxima of empirical processes, as given
by~\cite{klein2002inegalite} and~\cite{klein2005concentration}.

Our theory allows us to treat a number of new examples in which
concentration holds.  However, as (asymptotically) exact values for
the expectation of maxima of the empirical process are generally not
available, we cannot always provide explicit expressions for the point
of concentration in terms of the parameters of the model.

\subsection*{Set-up and notation}

Consider independent observations $X_1, \ldots , X_n$ taking values in
a space ${\cal X}$, a given class $\calF$ of real-valued functions on
${\cal X}$ and a non-negative regularization penalty $\pen : \ \calF
\rightarrow [0, \infty)$.  The empirical measure $P_n$ of a function
  $f : {\cal X} \rightarrow \R$ is defined by the average $P_n f \defn
  \frac{1}{n} \sum_{i=1}^n f(X_i)$, whereas the theoretical or
  population measure is given by \mbox{$Pf \defn \EE P_n f $.}

We let $\calF$ denote a class of loss functions, say indexed by a
parameter $g$ in a parameter space ${\cal G}$.  As a concrete example,
in the case of least-squares regression, the observations consist of
covariates along with real-valued responses of the form $\{(X_i,Y_i)
\}_{i=1}^n$.  In the least-squares case, the loss class takes the form
$\calF =\{ f_g(x,y) = (y - g(x))^2: \ g \in \calG \}$, where $\calG$
is some underlying collection of regression functions.

With this set-up, the regularized empirical risk estimator is
defined\footnote{In order to avoid digressions, we assume throughout
  this paper that ``argmin's''and ``argmax's'' exist and are unique in
  a suitable sense.} as
\begin{align}
\label{EqnERM}
 \hat f = \arg \min_{f \in \calF} \biggl \{ P_n f + \pen (f) \biggr
 \}.
\end{align}

We define the associated target function $f^0 \defn \arg \min
\limits_{f \in \calF} Pf$, corresponding to the population
minimizer, and we let
\begin{align}
\label{EqnDefnExcessRisk}
\tau^2(f) \defn P(f-f^0) + \pen (f)
\end{align}
the penalized excess risk.
  
In order to simplify the exposition, we give asymptotic statements at
places, using the classical scaling in which the sample size $n$ tends
to infinity.  For a sequence of positive numbers
$\{z_n\}_{n=1}^\infty$, we write
\begin{align*}
z_n = {\mathcal O} (1) \ {\rm if } \ \limsup z_n < \infty , \quad
\mbox{and} \quad z_n=o(1) \ {\rm if} \ z_n \rightarrow 0,
\end{align*}
as well as $z_n \asymp 1 $ if both $z_n = {\mathcal O} (1)$ and $1/z_n
= {\mathcal O}(1)$.  For two positive sequences $\{y_n\}_{n=1}^\infty$
and $\{z_n\}_{n=1}^\infty$, we write $z_n = {\mathcal O} (y_n)$ if
$z_n/y_n = {\mathcal O} (1)$, along with analogous definitions for the
$o$- and $\asymp$-notation.  We furthermore use the stochastic order
symbols are ${\mathcal O}_{\PP}$ and $o_{\PP}$.  In all our uses of
these forms of order notation, the arguments depend on $n$, but we
often omit this dependence to simplify notation.
  
With this set-up, the main results of this paper involve showing that
under certain conditions, we have
\begin{align*}
\biggl | { \tau ( \hat f ) - s_0 } \biggr | = o_{\PP} (s_0) . 
\end{align*}
Here $s_0$ is a deterministic quantity defined by the problem under
consideration; see equation~\eqref{s0.equation} for its precise
definition.  In our context, the result requires the complexity of the
problem to be in the nonparametric regime in which $\sqrt { \log n /n}
= o(s_0)$.  When a certain concavity condition is met, the $\log
n$-term can be removed.  This concavity condition holds in the normal
sequence model, as well as in all the examples given in Section
\ref{pure.section}.  In Section~\ref{direct.section}, there is no
$\log n$-term as well, but the concentration result there is for
$\sqrt{ P(\hat f - f^0)}$ as opposed to $\tau(\hat f)$.
  

\subsection*{Organization}
  
The remainder of this paper is organized as follows. In
Section~\ref{direct.section}, we provide a concentration result for
the normal sequence model and least squares with convex penalty, based
on a ``direct'' argument.  We then consider more general models and
loss functions, using the more indirect route originally taken by
\cite{chatterjee2014new}.  In
Section~\ref{theoretical-version.section}, we discuss the
deterministic counterpart of empirical risk minimization,
corresponding the population-level optimization problem.  Our theory
requires a certain amount of curvature of the objective function
around its minimum, a requirement that we term a second order margin
condition.  In Section~\ref{ERM.section}, we present a concentration
result (Theorem~\ref{local.theorem}) for a general loss function.
Section~\ref{quadratic.section} is devoted to a more careful analysis
of quadratic second-order margin conditions.
Section~\ref{pure.section} is devoted to the detailed analysis of two
examples in which the empirical process is linear in its
parameter---projection estimators for densities and linearized least
squares---whereas Section~\ref{exponential-families.section} provides
results for nonparametric estimation involving exponential families.
In Section~\ref{concentration.section}, we present the concentration
inequalities that underlie the proof of our indirect approach.  In
Section~\ref{shifted.section}, we provide a similar result as in
Section~\ref{ERM.section} but now for a shifted version of $\tau^2
(\hat f)$.  Finally, all proofs are in provided in
Section~\ref{proofs.section}.


\section{Direct approach to normal sequence model} 
\label{direct.section}

In this section, we analyze the concentration properties of
regularized least-squares estimators in the normal sequence setting.
The main contribution of this section is to provide a direct argument
that generalizes and sharpens the previous result
of~\cite{chatterjee2014new}.

Let $Y_i \in \R$ be a response variable and $X_i $ be a fixed
co-variable in some space ${\cal X}$, $i=1 , \ldots , n$.  The normal
sequence model is given by
\begin{align}
Y_i = g^0(X_i) + \epsilon_i \qquad \mbox{for $i=1 , \ldots , n$,}
\end{align}
where $\epsilon_1 , \ldots , \epsilon_n$ are i.i.d.\ mean-zero
Gaussians with variance $\sigma^2$, and the regression vector
\mbox{$g^0 \defn (g^0(X_1) , \ldots , g^0(X_n) )^T$} is unknown.  Let
us write the vector of responses as \mbox{$Y = (Y_1 , \ldots , Y_n
  )^T$,} and the noise vector as \mbox{$\epsilon \defn ( \epsilon_1 ,
  \ldots , \epsilon_n )$.}

Let $\pen : \R^n \rightarrow \R \cup \{ +\infty \} $ be a complexity
penalty, assumed to be convex.  The regularized least squares
estimator is given by
\begin{align}
\label{EqnRegularizedNS}
\hat g \defn \arg \min_{g \in \R^n } \biggl \{ \| Y - g \|_n^2 + {\rm
  pen} (g) \biggr \},
\end{align}
where for any vector $v \in \R^n$, we use the standard notation
\mbox{$\| v \|_n^2 \defn \frac{v^T v}{2} \; = \;
  \frac{\|v\|_2^2}{n}$.}

In past work, \cite{chatterjee2014new} analyzed the concentration of
the constrained variant of this estimator, given by
\begin{align}
\label{EqnConstrainedNS}
\hat g \defn \arg \min_{g \in \GSET} \biggl \{ \| Y - g \|_n^2 \biggr
\},
\end{align}
where $\GSET \subseteq \R^n$ is a closed, convex set.  Note that this
constrained estimator~\eqref{EqnConstrainedNS} is a special case of
the regularized estimator~\eqref{EqnRegularizedNS}, in which the
penalty function takes the form
\begin{align}
\label{EqnIndicator}
\pen(g) & \defn \begin{cases} 0 & \mbox{if $g \in \GSET$} \\ +\infty &
  \mbox{otherwise.}
\end{cases}
\end{align}

The following result guarantees that with for any convex penalty, the
estimation error $\| \hat g - g^0 \|_n$ of the regularized
estimator~\eqref{EqnRegularizedNS} is sharply concentrated around its
expectation $m_0 \defn \EE \| \hat g - g^0 \|_n$.
  
\begin{theorem}
\label{direct.theorem} 
For any convex penalty $\pen: \R^n \rightarrow \R \cup \{+\infty \}$,
the error in the regularized estimator~\eqref{EqnRegularizedNS}
satisfies
\begin{align}
 \PP \biggl ( \biggl | \| \hat g - g^0 \|_n - m_0 \biggr | \ge \sigma
 \sqrt {2t/n} \biggr ) \le \exp[-t] \qquad \mbox{for all $t > 0$.}
\end{align}
\end{theorem} 
See Section~\ref{proofs.section.direct} for the proof of this claim.
The argument is direct, using some basic facts from convex analysis,
and the concentration of Lipschitz functions of Gaussian vectors
(see~\cite{borell1975brunn}).
  
{\bf Remarks:} In terms of asymptotic behaviour, it follows that when
$\sigma = {\mathcal O} (1)$ and $1/ \sqrt {n} = o(m_0) $---the latter
condition corresponding to the non-parametric regime---it holds that
\begin{align*}
 \biggl | \| \hat g - g^0 \|_n - m_0  \biggr |   = o_{\PP} (m_0) .
\end{align*}
Moreover, it follows from its proof that Theorem~\ref{direct.theorem}
remains true if the population minimizer $g^0$ is replaced by any
other vector $g\in \R^n$. Thus, for instance, we can take $g$ as the
minimizer of the penalized noiseless problem
\begin{align*}
g^* \defn \arg \min_{g \in \R^n} \biggl \{ \| g - g^0 \|_n^2 + \pen
(g) \biggr \}
\end{align*}
With this choice, we also have concentration of $\| \hat g - g^* \|_n$
around its expectation $\EE \| \hat g - g^* \|_n $.

With the choice~\eqref{EqnIndicator}, the result also applies to the
constrained least squares estimate~\eqref{EqnConstrainedNS}, for any
closed convex set $\GSET \subseteq \R^n$.  In this context,
Theorem~\ref{direct.theorem} sharpens the previous result
of~\cite{chatterjee2014new}.

  
\section{Theoretical version of  minimization problem}
\label{theoretical-version.section}
  
For the remainder of the paper, we study the general empirical risk
minimizer.  Let $\tau (\hat f)$ be the excess
risk~\eqref{EqnDefnExcessRisk} associated with the empirical minimizer
$\hat f$ from equation~\eqref{EqnERM}.  We define the minimum possible
excess risk
\begin{align*}
\tau_{\rm min}^2 \defn \min_{f \in \calF} \tau^2 (f) = \min_{f \in
  \calF} \biggl \{ P( f-f^0) + \pen (f) \biggr \} .
\end{align*}
For all $s \ge \tau_{\rm min}$, we define the set $\calF_s \defn \{ f
\in \calF \, \mid \, \ \tau(f) \le s \}$, and the functions
\begin{align*}
\hat {\bf E}_n (s) = \max_{f \in \calF_s } ( P_n - P) (f^0-f), \quad
\mbox{and} \quad \ {\bf E}(s) \defn \EE( \hat {\bf E}_n (s) ) .
\end{align*}
In addition, we define the minimizers
\begin{align}
\label{s0.equation}
\hat s \defn \arg \min_{s\ge \tau_{\rm min}} \{ s^2 - \hat {\bf E}_n
(s)\} , \quad \mbox{and} \quad s_0 \defn \arg \min_{s\ge \tau_{\rm
    min} } \{ s^2 - {\bf E}(s) \} .
\end{align}
From hereon, we refer to minimizing the function $s \mapsto s^2 - {\bf
  E}(s)$ as the theoretical problem, and to minimizing $s \mapsto s^2
- \hat {\bf E}_n (s)$ as the empirical problem.  That the latter
indeed yields the risk associated with the original ERM
estimate~\eqref{EqnERM} is guaranteed by the following lemma.
 
\begin{lemma} \label{minimum.lemma}
We have\footnote{Recall that throughout this paper, we tacitly assume
  that the minimizers are unique.}  $ \tau(\hat f) = \hat s .$
  \end{lemma}
See Section~\ref{proof.minimum.lemma} for the proof of this claim.

In order to prove concentration results, we need a certain amount of
curvature of the function $s \mapsto s^2 - {\bf E}(s)$ around its
minimum.  The following condition is known as a second-order margin
condition; we also introduce a first-order margin condition in
Definition~\ref{1order-margin.definition}.
 
\begin{definition} \label{2order-margin.definition} 
Let $G$ be a strictly convex increasing function with $G(0) =0$ and
let $\underline \delta\ge 0$ and $\tau_{\rm max} \ge s_0 + \underline
\delta$ be constants.  If $s_0 > \tau_{\rm min} $ we say that the
second order margin condition holds in the range $ [ \tau_{\rm min} ,
  s_0 ) \cup ( s_0 + \underline \delta, \tau_{\rm max} ]$ with margin
function $G$ if for all $s \in [ \tau_{\rm min} , s_0 ) \cup ( s_0 +
  \underline \delta, \tau_{\rm max} ]$
 \begin{equation}
\label{2order-margin.equation}
  [s^2 - {\bf E}(s) ]- [s_0^2 - {\bf E}(s_0) ] \ge G( |s- s_0|) .
  \end{equation}
 If inequality~\eqref{2order-margin.equation} is only true for all $s
 \in ( s_0+ \underline \delta , \tau_{\rm max} ] $, then we say that
   the right-sided second order margin condition holds in the range
   $(s_0 + \underline \delta, \tau_{\rm max} ]$ with margin function
     $G$.  If no reference is made to any range, it means that the
     condition holds $\underline \delta =0$ and $\tau_{\rm max} =
     \infty$.
    \end{definition} 
    
 In the two-sided case in Definition \ref{2order-margin.definition},
 we do not allow for a gap for values of $s$ to the left of
 $s_0$. This choice is merely for simplicity and corresponds to our
 examples.

An important special case is quadratic margin behaviour, as formalized
in the following:
  \begin{definition}  
\label{2order-quadratic-margin.definition} 
If $s_0 > \tau_{\rm min}$ we say that the second order quadratic
margin condition is met in the range $[\tau_{\rm min} , s_0) \cup (
  s_0 + \underline \delta , \tau_{\rm max} ] $ with margin constant
$c>0$ if for all $s \in [\tau_{\rm min} , s_0) \cup ( s_0 + \underline
  \delta , \tau_{\rm max} ] $, it holds that
 \begin{align*}
\biggl [ s^2 - {\bf E}(s) \biggr ] - \biggl [ s_0^2 - {\bf E}(s_0)
  \biggr ] \ge (s-s_0)^2 / c^2 .
 \end{align*}
If this is only true for $s \in (s_0 + \underline \delta , \tau_{\rm
  max} ]$, then we say that the right-sided second order quadratic
  margin condition holds in the range $(s_0 + \underline \delta ,
  \tau_{\rm max}]$ with margin constant $c>0$.
\end{definition} 
  
When the two-sided condition in
Definition~\ref{2order-quadratic-margin.definition} holds with
$\underline \delta =0$, then it corresponds to a form of \emph{strong
  convexity} of the function $s \mapsto s^2 - {\bf E}(s)$ at $s_0$.
  
Clearly, if $s \mapsto {\bf E}(s)$ is concave then $s \mapsto s^2 -
{\bf E}(s)$ is strictly convex.  In fact, then the second order
quadratic margin condition holds with margin constant $c$ at most
equal to $1$.  This type of condition holds in the normal sequence
setting, as exploited by~\cite{chatterjee2014new}.  In the latter
paper, the map $s \mapsto \hat {\bf E}_n (s)$ is concave and hence
then also $s \mapsto {\bf E}(s)$ is concave. Moreover, the empirical
function $ s \mapsto s^2 - \hat {\bf E}_n (s)$ is then convex, which
allows one to remove the $\log n$-factor.  We will consider conditions
for (right-sided) second order quadratic margin behaviour in
Section~\ref{quadratic.section}.
  

\section{Concentration of ERM}
\label{ERM.section}
 
We now turn to the statement of our main result on concentration of
ERM in the general setting.  We begin by specifying some conditions
that underlie the result.  First, we require a uniform boundedness
condition:
\begin{condition} 
\label{uniformly-bounded.equation}
The function class $\calF$ is uniformly bounded, meaning that
\begin{align*}
K \defn \max_{f \in \calF } \| f - f^0 \|_{\infty} < \infty .
\end{align*}
\end{condition}
We note that this condition can be removed if one first shows that,
for a suitable constant $K$, the minimizer $\hat f$ satisfies the
bound $\| \hat f- f^0 \|_{\infty} \le K$ with high probability.
  
When Condition~\ref{uniformly-bounded.equation} holds, one may take
\begin{align*}
\tau_{\rm max}^2 \defn 2 K + \pen (f^0) .
\end{align*}
However, in order to obtain a sharper result, one may first want to
prove that $\tau^2(\hat f)$ is much smaller than \mbox{$2 K + \pen
  (f^0)$} with high probability.  In fact, there is a substantial
literature on techniques for showing that $ \tau(\hat f) = {\mathcal
  O}_{\PP} (s_0)$ for various problems.  As we discuss in
Section~\ref{shifted.section}, similar results exist for the shifted
version, in particular showing that $\tau^2(\hat f) - \tau_*^2 =
{\mathcal O}_{\PP} (s_0^2 - \tau_*^2 )$, where $\tau_*^2$ is a
suitably chosen constant.

We now come to our ``first order'' margin condition, which quantifies
the curvature of $f \mapsto Pf$ around its minimum.  To avoid confusion
with the ``second order" margin condition we call it a ``curvature condition".

%
%

For $f \in \calF$, define the
variance
\begin{align*}
\sigma^2 (f) \defn {1 \over n} \sum_{i=1}^n \biggl [ \EE f^2 (X_i) -
  (\EE f(X_i))^2 \biggr ] . 
\end{align*}

\begin{definition} 
\label{1order-margin.definition}
A quadratic curvature condition with constant $C>0$ is said to
hold if 
\begin{align}
\label{EqnQuadMargin}
 P(f-f^0) \ge \frac{\sigma^2(f-f^0)}{C^2} \quad \mbox{for all $f \in
   \calF$.}
\end{align}
\end{definition}
  
We take the quadratic curvature condition as basis for our results. An
extension to more general curvature is omitted here to avoid
digressions.
    
Let ${\cal J}$ be a strictly increasing function defined on
$[\tau_{\rm min} , \infty)$ and such that ${\cal J} (\tau_{\rm min})
  =0 $.  We then define a new function
\begin{align}
\label{EqnDefnPhi}
\Phi_{\cal J} (u) \defn [ {\cal J}^{-1} (u)]^2 \qquad \mbox{for all $u
  >0$.}
\end{align}
Our choice of the square here is linked to the quadratic curvature
condition~\eqref{EqnQuadMargin}.

\begin{condition} \label{r0.condition}
There is a constant $m_n$, a strictly increasing function $\cal J$
such that the function $\Phi_{\cal J}$ is strictly convex, and such
that the bound
\begin{align}
{\bf E}(s) \le \frac{{\cal J} (s)}{m_n}
\end{align}
holds for all $s \ge \tau_{\rm min}$.
\end{condition}

We also define the convex conjugate function $\PhiConj_{\cal J}(v)
\defn \sup \limits_{u > 0} \Big \{ v u - \Phi_{\cal J}(u) \Big \}$,
and make use of the Fenchel-Young inequality
\begin{align*}
 uv \le \Phi_{\cal J}(u) + \PhiConj_{\cal J} (v) , \qquad \mbox{valid
   for all pairs $u,v >0$.}
\end{align*}
Finally, in terms of the previously defined quantities, we define
\begin{align*}
 r_0^2 \defn 2 C^2 \PhiConj_{\cal J} \big( 4 K / ( m_n C^2 ) \big).
\end{align*}
With this notation, the following theorem is our main result:
\begin{theorem}
\label{local.theorem} 
Suppose 
that:
\bcar
\item Conditions~\ref{uniformly-bounded.equation}
  and~\ref{r0.condition}, as well as the quadratic curvature condition
  with constant $C$ hold.
\item The right-sided second order
margin condition holds in the range $ ( s_0 + \underline \delta,
\tau_{\rm max} ] $ with margin function $G$, with associated convex
conjugate $\Gconj$. 
\ecar
Then there is a constant \mbox{$c_0 = c_0(C,
  K)$,} and a function $\delta:(0, \infty) \rightarrow (0, \infty)$
such that
\begin{multline*}
G(\delta(t))  \leq \Gconj \biggl ( c_0 \sqrt {[ t + \log (1+ \sqrt
    {n\tau_{\rm max}^2} )]/n} ) \biggr ) \\ + c_0 \biggl ( (s_0+r_0)
\sqrt {[ t + \log (1+ \sqrt {n\tau_{\rm max}^2} )]/n} + [ t + \log (1+
  \sqrt {n\tau_{\rm max}^2} )]/n \biggr ),
\end{multline*}
and such that the following deviation inequality holds:
\begin{align*}
 \PP \biggl (\tau_{\rm max} \ge \tau(\hat f) > s_0 + \max\{ \underline
 \delta, \delta (t)\} \biggr ) \le \exp[-t] \qquad \mbox{for all $t >
   0$.}
\end{align*}
If, in fact, the two-sided version of the second order margin
condition holds over $[\tau_{\rm min} , s_0 ) \cup (s_0 + \underline
  \delta , \tau_{\rm max} ]$, then we have
\begin{align*}
 \PP \biggl ( |\tau(\hat f) - s_0 |> \max\{ \underline \delta, \delta
 (t) \} , \ \tau(\hat f) \le \tau_{\rm max} \biggr ) \le 2 \exp[-t]
 \qquad \mbox{for all $t > 0$.}
\end{align*}
\end{theorem}

{\bf Asymptotics:} If the second order margin condition holds with the quadratic
function $G(u)= \frac{u^2}{2 c^2}$, then its convex conjugate
$\Gconj(v) = \frac{c^2 v^2}{2} $ is also quadratic.  Thus, under the
scalings $C = {\mathcal O} (1)$, $K = {\mathcal O} (1)$ and $r_0
\asymp s_0$, we then find that
\begin{align*}
\delta (t) = {\mathcal O} \biggl ( ({\log n / n} )^{1/2} + ({ s_0^2
  \log n / n}) ^{1/4} \biggr )
\end{align*}
for each fixed $t$.  Hence, whenever $\sqrt {\frac{\log n}{n}} =
o_{\PP} (s_0)$, then we are guaranteed that \mbox{$\delta (t) =
  o_{\PP} (s_0) $.}


\section{Second order quadratic margin behaviour}
\label{quadratic.section}

In this section, we investigate conditions under which the
(right-sided) second order quadratic margin condition holds over an
appropriate range.  In particular, we extend the setting of
\cite{chatterjee2014new} to the case where one has a strictly convex
penalty in Lemma \ref{concaveq.lemma}, and to approximate forms of
concavity in Lemmas~\ref{approx-concave.lemma}
and~\ref{concave.lemma}. We note that it is possible to formulate
different results with other combinations of conditions, but we
omit this here.
\begin{lemma}
\label{concaveq.lemma}
Let $\calF \defn \{ f_g \mid \ g \in \calG \}$ be a class of loss
functions indexed by a parameter $g$ in a parameter space
$\calG$. Assume $\calG$ is a convex subset of a linear vector space,
the mapping $ g \mapsto f_g - Pf_g$ is linear, and that for some $q >
1$, the mapping $g \mapsto \tau^{2/q} (f_g)$ is convex.  For some
constant $M > 0$, define \mbox{$\tau_{\rm max} \defn (M+1) s_0$.}
Then the right-sided second order quadratic margin condition holds in
the range $(s_0 , \tau_{\rm max} ] $ with constant
\begin{align*}
c= \sqrt { 2q^{-1} (q-1) (M+1)^{-{2(2-q) \over q} } }.
\end{align*}
Moreover, when $q=2$ and $s_0 >\tau_{\rm min} $, then the (two-sided)
second order quadratic margin condition holds with $c=1$.
  \end{lemma}
We note that the latter two-sided second order quadratic margin
condition corresponds to the favourable setting of the normal sequence
model, as studied by \cite{chatterjee2014new}.

{\bf Asymptotics:} The idea in the above lemma is that one first
proves by separate means that $\tau(\hat f) = {\mathcal O}_{\PP}
(s_0)$. There is a large literature on bounds of this form; for
example, see~\cite{koltchinskii2011oracle} and references therein.
One can then take $M= {\mathcal O} (1)$.

We sometimes write $\hat {\bf E}_n(\cdot) =: \hat {\bf E}_n^{\tau}
(\cdot)$ and ${\bf E}(\cdot) =: {\bf E}^{\tau} (\cdot)$ so as to
highlight their dependence on $\tau$.  For $f \in \calF$, define the
functionals
\begin{align*} 
\varsigma^2 (f) \defn {\rm c}^2 \sigma^2 (f-f^0)+ { \pen (f)}, \quad
\mbox{and} \quad \varsigma_{\min} \defn \min_{f \in \calF} \varsigma
(f),
\end{align*}
where ${\rm c} > 0 $ is some constant.  Moreover, for $s \ge
\varsigma_{\rm min}$, let us define
\begin{align*}
\hat {\bf E}_{n}^{ \varsigma} (s) \defn \max_{f \in \calF: \ \varsigma
  (f) \le s} (P_n - P) (f^0 - f) , \quad \mbox{and} \quad {\bf
  E}^{\varsigma} (s) \defn \EE \hat {\bf E}_{n}^{ \varsigma} (s).
\end{align*}
  
\begin{lemma} 
\label{approx-concave.lemma}
Suppose that the function $s \mapsto {\bf E}^{\varsigma}(s)$ is
concave, and that
\begin{align*}
\frac{\varsigma^2 (f)}{A} \; \le \; \tau^2 (f) \; \le \; A
\varsigma^2(f) , \quad \mbox{for all $f \in \calF$, and $\tau_{f} \le
  \tau_{\max}$,}
\end{align*}
where $A^2 = 1+ \epsilon $ for some $ \epsilon >0$ satisfying $\sqrt
\epsilon (1+ \epsilon) < 1/2 $. Let \mbox{$\tau_{\rm max} \defn (M+1)
  s_0$} for some $M > 0$ and
\begin{align*}
\underline \delta \defn 2 \big [ \sqrt \epsilon ( 2 \sqrt \epsilon M +
  1) \big ]^{1/2} s_0 .
\end{align*}
Then when $s_0 > \tau_{\rm min}$, the quadratic second order margin
condition holds in the range $[ \tau_{\rm min} , s_0) \cup ( s_0 +
  \underline \delta, \tau_{\rm max} ] $ with constant $c=4$.
  \end{lemma} 
  
  {\bf Asymptotics} As in Lemma \ref{concaveq.lemma} one may first
  prove by separate means that $\tau(\hat f) = {\mathcal O}_{\PP}
  (s_0)$ and then take $M= {\mathcal O} (1)$.

Lemma~\ref{approx-concave.lemma} requires the function ${\bf
  E}^{\varsigma}$ to be concave.  We now present conditions under
which this is indeed the case.

\begin{lemma} \label{concave.lemma} 
Let $\calF \defn \{ f_g \mid \ g \in \calG \}$ be a class of loss
functions indexed by the parameter $g$ in a parameter space
$\calG$. Assume $\calG$ is a convex subset of a linear vector space,
and that $g \mapsto \sqrt {\pen (f_g)} $ is convex and $ g \mapsto f_g
- Pf_g$ is linear.  Then the function $s \mapsto \hat {\bf E}_{n}^{
  \varsigma} (s)$ is concave.
\end{lemma}

In fact, we show concavity of the empirical version $\hat {\bf
  E}_n^{\varsigma} $, which then implies concavity of ${\bf
  E}^\varsigma_n$.The reasoning is along the lines
of~\cite{chatterjee2014new}, and of the corresponding part of the
proof of Lemma \ref{concaveq.lemma}.


\section{Some ``pure'' cases}
\label{pure.section}

In this section, we examine a number of problems that are ``pure'' in
the sense that the empirical process enters in a linear manner.  The
simplest example of such a pure case is the normal sequence model
studied in Section~\ref{direct.section}, and we examine some other
examples here.

More precisely, consider a class of the form $\calF \defn \{ f_g \mid
\ g \in \calG \}$, where $\calG$ is a convex subset of a normed linear
vector space $( \bar\calG , \| \cdot \| ) $.  The pure case
corresponds to problems in which the mapping $g \mapsto f_g - Pf_g$ is
linear, and moreover, we have $P(f_g-f^0 ) = \| g - g^0 \|^2 $, where
$g_0 = \arg \min \limits_{g \in \calG} P f_g$, which ensures the
equivalence $f^0= f_{g^0}$.


\subsection{Density estimation using projection}
\label{projection.section}

Let $X_1 , \ldots , X_n$ be i.i.d.\ random variables with distribution
$P$ taking values in a space ${\cal X}$.  For a sigma-finite measure
$\nu$ on ${\cal X}$, let $\| \cdot \|$ denote the $L^2(\nu)$-norm.
Let $\calG$ be a convex subset of a linear vector space ${\bar G}
\subset L^2(\nu) $, and suppose that density $g^0 \defn dP/ d \nu $ is
a member of the model class $ \calG $.  With this set-up, we consider
the estimator
\begin{align*}
\hat g \defn \arg \min \biggl \{ - P_n g + \frac{\| g \|^2}{2} +
\lambda^2 I^q (g) \biggr \},
\end{align*}
where $I$ denotes some pseudo-norm on ${\bar G}$, the exponent $q \in
(1,2]$, and $\lambda \geq 0$ is a regularization parameter.

In order to analyze the concentration properties of this estimator
using our general theory, we begin by casting it within our framework.
For each $g \in \calG$, define $f_g \defn -g+ \frac{1}{2} \| g\|^2 $,
as well as the associated function class \mbox{$\calF \defn \{ f_g
  \mid \ g \in \calG \} $.}  With these choices, for all $g \in
\calG$, we have $(P_n-P)f_g = -(P_n-P)g$, and moreover
\begin{align*}
P( f_g - f_{g_0}) = -P(g-g^0) + \frac{\| g \|^2}{2} - \frac{\| g^0
  \|^2}{2} \, = \, \| g - g^0 \|^2 .
\end{align*}

We split our analysis into several cases, depending on the nature of
the \mbox{penalty $I$.}

\subsubsection{ Case 1: No penalty}

In this case, we assume:\\
$\circ$ Condition \ref{uniformly-bounded.equation} holds with $K=
{\mathcal O} (1)$,\\ 
$\circ$ ${\bf E}(s) \le {\cal J} (s) /\sqrt n $,
\\
$\circ$ ${\cal J} (s) = A s^{1- \alpha} $, $\exists \ A = {\mathcal O}
(1) $, $\exists \ 0< \alpha < 1 $ not depending on $n$,\\
$\circ$ $s_0 \asymp n^{-{1 \over 2 (1+ \alpha)}} $.

It then follows from Theorem~\ref{local.theorem} combined with
Lemma~\ref{concaveq.lemma} that
\begin{align*}
\biggl |{ \| \hat g - g^0 \| - s_0 \over s_0 } \biggr | = {\mathcal
  O}_{\PP}\biggl ( \sqrt {\log n \over n^{\alpha \over 1+ \alpha}
}\biggr ) = o_{\PP} (1) .
\end{align*}
In fact, in this case, the $\log n$-term can be removed because
Lemma~\ref{concave.lemma} ensures that the map $s \mapsto \hat {\bf
  E}_n (s)$ is concave.


\subsubsection{ Case 2: Quadratic penalty}

In this case, we assume: \\
$\circ$ $q=2$, \\ 
$\circ$ Condition \ref{uniformly-bounded.equation} with $K= {\mathcal
  O} (1)$,\\
$\circ$ ${\bf E}(s) \le {\cal J} (s)/m_n $, \\
$\circ$ ${\cal J} (s) = A s$, $m_n = \sqrt n \lambda^{\alpha} $,
$\exists \ A = {\mathcal O} (1) $, $\exists \ 0< \alpha < 1 $ not
depending on $n$,\\
$\circ$ $s_0 \asymp 1 / (\sqrt n \lambda^{\alpha} ) $,\\
$\circ$ $ \lambda^{2\alpha} \log n = o(1) $.

Then Theorem \ref{local.theorem} combined with
Lemma~\ref{concaveq.lemma} implies that
\begin{align*}
\biggl | { \tau (\hat f) - s_0 \over s_0 } \biggr | = {\mathcal
  O}_{\PP} \biggl ( \sqrt {\lambda^{2 \alpha} \log n } \biggr ) =
o_{\PP} (1) . 
\end{align*}
As before, the $\log n $-factor can be removed.

One sees from the condition $\lambda^{2 \alpha } \log n = o(1)$ that
the regularization parameter $\lambda$ must be sufficiently small.
Moreover, hiding in our conditions is the fact that the penalty
$I(f^0)$ is not too large.  Indeed, we have $s_0^2 \asymp {1 \over n
  \lambda^{2 \alpha } } \ge \tau_{\rm min}^2$.  Consequently, when
$\tau_{\rm min}^2 \asymp \pen (f^0) = \lambda^2 I^2 (f^0)$, we must
have
\begin{align}
\label{bound-I.equation} 
I^2(f^0) = { \mathcal O} \Big ( \big ( { n \lambda^{2(1+ \alpha) } }
\big )^{-1} \Big) .
\end{align}
As a special case, suppose that we take $\lambda = n^{-{1 \over 2(1+
    \alpha)}} $.  Then we find \mbox{$s_0 \sim n^{-{1 \over 2( 1+
      \alpha) }}$} as in the previous section.  Then the
bound~\eqref{bound-I.equation} yields $I^2(f^0) ={ \mathcal O} (1) $.
Otherwise, we see that $\tau (\hat f) $ concentrates on the boundary
$\tau_{\rm min}$, and that in this example, we have
\begin{align*}
\tau_{\rm min}^2 = \min_{g \in \calG } \biggl \{ \| g - g^0 \|^2 +
\lambda^2 I^2 (g) \biggr \} .
\end{align*}


\subsubsection{Case 3: Strictly convex penalty}

In this case, we assume:\\ 
$\circ$ $1<q\le 2 $ not depending on $n$,\\ 
$\circ$ Condition \ref{uniformly-bounded.equation} with $K= {\mathcal
  O} (1)$,\\ 
$\circ$ ${\bf E}(s) \le {\cal J} (s) / m_n $, \\ 
$\circ$ ${\cal J} (s) = A s^{1+ (2/q-1) \alpha} $, $m_n= \sqrt n
\lambda^{2\alpha/q} $, $\exists \ A = {\mathcal O} (1) $, $\exists
\ 0< \alpha < 1 $ not depending\\ \ \ \ on $n$,\\
 $\circ$ $s_0 \asymp (\sqrt n \lambda^{2\alpha/q} )^{-{q \over q-
    (2-q) \alpha}} $,\\
$\circ$ $ \log n / ( n^{2-q} \lambda^4)^{\alpha \over q- (2-q)\alpha}
= o(1) $.

Under this condition, Theorem~\ref{local.theorem} combined with
Lemma~\ref{concaveq.lemma} implies the deviation result
\begin{align*}
{ \tau(\hat f) - s_0\over s_0 } \le z_n, \ {\rm where} \ z_n =
{\mathcal O}_{\PP} \biggl ( { \log n \over n s_0^2 } \biggr ) +
o_{\PP} (1) = o_{\PP} (1) .
\end{align*}


\subsection{Linearized least squares regression}
\label{linearizedLS.section}

Let $\{ (X_i , Y_i)\}_{i=1}^n $ be i.i.d. samples taking values in
$\R^p \times \R$.  We assume the model
\begin{align*}
Y_i = g^0 (X_i) + \epsilon_i , \quad \mbox{for $i=1 , \ldots , n$,}
\end{align*}
where $\epsilon_i \sim {\cal N} (0,1) $ is independent of $X_i$, and
the function $g^0$ belongs to a convex model class $\calG$. Assume
$\calG$ is a convex subset of a linear vector space $\bar \calG$, and
moreover that
\begin{align*}
K_X \defn \max_{g \in \calG } \| g - g^0 \|_{\infty} < \infty . 
\end{align*}
We moreover assume $K_0 \defn \| g^0 \|_{\infty} < \infty$.  Let $I$
be some pseudo-norm on $\bar \calG$, we consider the estimator
\begin{align}
\label{EqnLinearizedLeastSquares}
\hat g \defn \arg \min_{g \in \calG } \biggl \{ - \frac{1}{n}
\sum_{i=1}^n Y_i g(X_i) + \frac{P g^2}{2} + \lambda^2 I^2 (g) \biggr
\} .
\end{align}
Note that implementation of this estimator requires that $P g^2$ is
known or can be computed for all $g \in \calG$.

Given the form of the estimator~\eqref{EqnLinearizedLeastSquares},
we have
\begin{align*}
\calF = \big \{ f_g (x,y) = -y g(x) + Pg^2 /2 \, \mid \, \ g \in \calG
\big \}.
\end{align*}
This class has an envelope function $F$ satisfying the conditions of
Lemma~\ref{kleinrio-truncated.lemma}. To see this, note the bounds
\begin{align*}
\max_{g \in \calG} | \epsilon_i (g(X_i) - g^0 (X_i) | & \le
|\epsilon_i| K_X, \\ 
\max_{g \in \calG } |g^0 (X_i) (g(X_i) - g^0 (X_i)) | & \le K_X K_0,
\quad \mbox{and} \\
\max_{g \in \calG } \big | Y_i (g(X_i) - g^0 (X_0)) \big | & \le
(|\epsilon_i| + K_0 )K_X .
\end{align*}
 Lemma
 \ref{kleinrio-truncated.lemma} can now be used as concentration tool.
 Still, as the class $\calF$ is not uniformly bounded in this case
 one cannot apply Lemma \ref{Klein-Rio.lemma}. The strategy may then
 be to first prove that $\tau(\hat f) = {\mathcal O}_{\PP} (s_0)$ and
 then that ${\bf E}(s) \asymp s_0$ for $s \asymp s_0$.  The results
 are then as in the previous subsection albeit that we are facing an
 additional $\log n$ factor.
 
Linearized least-squares regression refers to the special case of a
linear model.  More concretely, define the design matrix $X \in \R^{n
  \times p}$ rows $X_i^T$ for $i=1, \ldots , n $, as well as the
sample covariance matrix $\hat \Sigma \defn X^T X / n $.  The linear
model consists of the function class
\begin{align*}
 \calG = \{ g_{\beta} (x) = x^T \beta :\ \beta \in {\cal B} \},
\end{align*}
along with the true function $g^0 (x) = x^T \beta^0$.  Our conditions
then require that the population covariance matrix $\Sigma_0 \defn \EE
\hat \Sigma $ is known, and moreover that ${\cal B}$ is a convex
subset of $\R^p$ and satisfying, for some constants $K_0$ and $K_X$,
\mbox{$| X^T_1 \beta^0 | \le K_0$} and \mbox{$|X^T_1 (\beta - \beta^0
  ) | \le K_X$} for all $\beta \in {\cal B} $.  The latter is for
example true when $\| X_1\| \le K_X$ and $\| \beta - \beta^0 \|_1 \le
1 $ (say) for all $\beta \in {\cal B}$.


\section{Exponential families with squared norm penalty}
\label{exponential-families.section}

We now turn to some examples involving exponential families.
Throughout this section, we specialize to the case of squared norm
penalties, noting that more general penalties can be studied as in the
previous section.

\subsection{Density estimation}
 
Suppose that $X_1 , \ldots , X_n$ are i.i.d. random variables with
distribution $P$ taking values in ${\cal X}$.  Given a sigma-finite
measure $\nu$ on ${\cal X}$, let us define the function class
\begin{align*}
\bar \calG \defn \biggl \{ g: \ {\cal X} \rightarrow \R: \ \int \exp [
  g(x)] d \nu (x) < \infty \biggr \},
\end{align*}
and note that it is convex.  We define a functional on $\bar \calG$
via
\begin{align*}
\dfun(g) \defn \log \biggl ( \int \exp [g(x) ] d \nu(x) \biggr ) .
\end{align*}
Let $\calG \subset \bar \calG$ be a convex subset and define the
function $f_g = -g + \dfun (g)$, for each $g \in \calG$, along with
the associated function class $\calF \defn \{ f_g \mid g \in \calG
\}$.  Letting $I$ be a pseudo-norm on ${\bar G}$ along with with a
non-negative regularization weight $\lambda$, we consider the
estimator
\begin{align*}
 \hat g \defn \arg \min_{g \in \calG } \biggl \{ -P_n g + {\rm d} (g)
 + \lambda^2 I^2 (g) \biggr \},
\end{align*}
and define $\hat f \defn f_{\hat g}$.

For identification purposes, we take the functions in $\calG$ to be
centered--- that is, such that $ \int g d \nu =0$ for each $g \in
\calG$.  For simplicity and without loss of generality, we take $g^0
\equiv 0 $ so that $f^0 \equiv 1 $ and $\nu=P$.  Since $P$ is unknown,
the centering of the functions (in actual practice) is done with
respect to some other measure.  This difference does not alter the
theory but should be kept in mind when examining the assumptions.
 
The following lemma relates the function $\dfun$ to the second moment:
 \begin{lemma} \label{exponential-family.lemma} 
Suppose that $K\defn \max_{g \in \calG } \| g \|_{\infty} < \infty$.
Then we have
\begin{align*}
\max_{g \in \calG } \left | { \dfun (t g) \over t^2 P g^2 } -{1 \over
  2} \right | = {\mathcal O} (t ), \qquad \mbox{valid as $t \downarrow
  0$.}
\end{align*}
  \end{lemma}
As a useful corollary, it gives us an asymptotic expression for
$\dfun(g)$ as the \mbox{$\ell_\infty$-norm} of $g$ shrinks.  In
particular, we let define $\calG_{\infty} (\eta) \defn \{ g :\ \| g
\|_{\infty} \le \eta \}$.
\begin{corollary}
\label{CorEta}
For each $g \in \calG_\infty(\eta)$, we have
\begin{align*}
{\rm d} ( g) = {1 \over 2} P g^2 (1+ {\mathcal O} (\eta ) ), \qquad
\mbox{valid as $\eta \downarrow 0$.}
\end{align*}

\end{corollary}
 
We are now equipped to state a result.  Suppose that:\\
$\circ$ $r_0 = {\mathcal O} (s_0)$,\\ $\circ$ $ \tau (\hat f) =
{\mathcal O}_{\PP} ( s_0) $,\\ $\circ$ $ \| \hat g \|_{\infty} =
o_{\PP} (1) $,\\ $\circ$ $\sqrt { \log n / n } = o_{\PP} (s_0)
$. \\ 

Combining Theorem~\ref{local.theorem} with
Lemmas~\ref{approx-concave.lemma} and~\ref{concave.lemma}, along with
Corollary~\ref{CorEta}, guarantees that $\tau(\hat f) $ concentrates
on $s_0$, and in particular, we have
\begin{align*}
 | \tau (\hat f) - s_0 | = o_{\PP} (s_0) .
\end{align*}

 
\subsection{ Regression with fixed design}

Let $\{ (X_i , Y_i) \}_{i=1}^n$ be independent observations taking
values in the Cartesian product space ${\cal X} \times \R$.  We assume
the design $\{ X_i \}_{i=1}^n$ is fixed, and that for some given
sigma-finite measure $\nu$, the log-density of $Y_i$ given $X_i$ takes
the log-linear form
\begin{align*}
 -Y_i g^0(X_i) + \dfun ( g^0(X_i)),
\end{align*}
where the function $\dfun ( \xi) \defn \log \big ( \int \exp[ y \xi ]
d \nu \big )$ has domain 
\begin{align*}
\Xi \defn\{ \xi \in \R \, \mid \, \ \int \exp [ y \xi ] d \nu < \infty
\}.
\end{align*}

We define $g\defn( g(X_1 ) , \ldots , g(X_n))^T \in \R^n $, let
$\calG$ be a convex subset of $\R^n$ and use, for $v \in \R^n$, the
notation $\| v \|_n^2 \defn v^T v / n $.  Letting $I$ be a pseudo-norm
on $\R^n$, we consider the estimator
 \begin{align*}
\hat g \defn \arg \max_{g \in \calG }\biggl \{ -Y^T g /n +
\sum_{i=1}^n \dfun (g (X_i) )/n + \lambda^2 I^2 (g) \biggr \} .
 \end{align*}
 In this case, the effective function class takes the form
\begin{align*}
\calF = \{ f_g (x, y) = -y g(x) + \dfun ( g(x)) : \ g \in \calG
\},
\end{align*}
and we have $ \hat f = f_{\hat g}$.
 
Let us assume that: \\ 
$\circ$ $ \sup_{\tau^2 (f_g) \le \eta }{ \sum_{i=1}^n \dfun ( g(X_i))
  - \dfun (g^0 (X_i) ) ) / (n \| g - g^0 \|_n^2) } = {\rm c}^2 (1+ o(
1 ) ) $, $\eta \downarrow 0 $,\\ 
$\circ$ $ \| \hat g - g^0\|_n = o_{\PP} (1)$,\\
$\circ$ $r_0 = {\mathcal O} (s_0)$,\\ 
$\circ$ $ \tau (\hat f) = {\mathcal O}_{\PP} ( s_0) $,\\ $\circ$
$\sqrt { \log n / n } = o_{\PP} (s_0) $. \\ 

Then Theorem~\ref{local.theorem} in conjunction with
Lemmas~\ref{approx-concave.lemma} and~\ref{concave.lemma} guarantees
that $\tau(\hat f) $ concentrates on $s_0$, and moreover that
\begin{align*}
 | \tau (\hat f) - s_0 | = o_{\PP} (s_0) .
\end{align*}


\section{Concentration for maxima of  empirical processes} 
\label{concentration.section}

The following result of \cite{klein2002inegalite} (see also
\cite{klein2005concentration}) is our main tool.
  
\begin{theorem} \label{Klein-Rio.theorem} Define
$K \defn \max_{f \in \calF_s } \| f - f^0 \|_{\infty}$ and $\sigma_s
  \defn \max_{f \in \calF_s } \sigma (f- f^0)$.  Then for all $t \ge
  0$, we have
\begin{align*}
\hat {\bf E}_n (s) & \ge {\bf E}(s) - \sqrt { 8 K {\bf E}(s) + 2
  \sigma_s^2 } \sqrt {t/n} - \frac{K t}{n}, \qquad \mbox{and} \\
\hat {\bf E}_n (s) & \le {\bf E}(s) + \sqrt { 8 K {\bf E}(s) + 2
  \sigma_s^2 } \sqrt {t/n} + \frac{2 K t}{3n},
\end{align*}
where each bound holds with probability at least $1- \exp[ -t]$,
\end{theorem}

We next present a consequence for the case where the functions in
$\calF_s$ are not uniformly bounded, but have a (sub-Gaussian)
envelope function. This result is invoked in the analysis of
Section~\ref{linearizedLS.section}.

\begin{lemma}\label{kleinrio-truncated.lemma} 
Assume that for some constants $c_F\ge 1$ and $C_F\ge 1$, the envelope
function $F (\cdot) \defn \max_{f \in \calF_s } | f(\cdot) - f^0
(\cdot ) |$ satisfies the bounds
\begin{equation}
P F^2 {\rm l } \{ F > t\}  \le c_F^2 \exp[-t^2 / C_F^2] , \ t >0 . 
 \end{equation}
Then for all $t>0$ with probability at least $1- \exp[-t] - 1/t^2 $
\begin{multline*}
\hat {\bf E}_n (s) \ge {\bf E} (s) - \sqrt { 8 C_F \sqrt {\log n} [
    {\bf E} (s) + 2 c_F (c_F +t) / n ] + 2 \sigma_s^2 } \sqrt {t/n}
\\- C_F t\sqrt {\log n } / n - c_F (4 c_F +t) / n,
\end{multline*}
and with probability at least $1- \exp[-t] - 1/t^2 $
\begin{multline*}
 \hat {\bf E}_n (s) \le {\bf E} (s) - \sqrt { 8 C_F \sqrt {\log n} [
     {\bf E} (s) + 2c_F^2 / n ] + 2 \sigma_s^2 } \sqrt {t/n} 
\\- 2C_F t \sqrt {\log n }  / (3n) -  c_F ( 4 c_F +t) / n .
\end{multline*}
\end{lemma}

%
 
In the next lemma, we replace the quantity ${\bf E} (s)$ appearing in
the square-root of Theorem~\ref{Klein-Rio.theorem} by a suitable upper
bound.
 
\begin{lemma} 
\label{Klein-Rio.lemma} 
Under Conditions~\ref{uniformly-bounded.equation}
and~\ref{r0.condition}, we have
\begin{align*}
 \hat {\bf E}_n (s) & \ge {\bf E}(s) - 2 C s \sqrt {t/n} + r_0 \sqrt
      {t/n} - Kt/n, \quad \mbox{and} \\
 \hat {\bf E}_n (s) & \le {\bf E}(s) + 2 C s \sqrt {t/n} + r_0 \sqrt
      {t/n} + 2Kt/(3n),
\end{align*}
where each bound holds with probability at least $1 - \exp[-t]$.
\end{lemma}

  
\section{The shifted version}
\label{shifted.section}

For a scalar $ \tau_*^2 \ge \tau_{\rm min}^2 $ to be chosen, we study
in this section the ``shifted'' function
\begin{align*}
 {\bf F}(s) \defn \max_{\tau^2 (f) \le \tau_*^2 + s^2 } (P_n - P) (
 f^0 - f), \qquad \mbox{defined for $s^2 \ge \tau_*^2 - \tau_{\rm
     min}^2$.}
\end{align*}
This shifted version may be of interest when $\tau^2(\hat f)$ is of
larger order than \mbox{$P(\hat f - f^0)$.} The idea is then to
replace $g^0$ in the previous sections by the function \mbox{$g^*
  \defn \arg \min_{g \in \calG }\tau^2 (g) $.} One then needs curvature
conditions on $R(g) - R(g^*)$ instead of $R(g) - R(g^0)$. This we
handle here by the notion of an ``oracle potential'', as defined in
Definition~\ref{potential.definition} below.


Lemma~\ref{shifted.lemma} shows that curvature conditions on the
function $\QUADE{s} \defn s^2 - {\bf E}(s)$ are weaker than those on
the function $\QUADF{s} \defn s^2 - {\bf F}(s)$.  Using the shorthand
$s_*^2 = s_0^2 - \tau_*^2$, the following lemma summarizes this fact:

\begin{lemma} \label{shifted.lemma} 
For any $s \ge \tau_{\rm min}$ and $\tilde s^2= s^2 - \tau_*^2$, we
have
\begin{align*}
\QUADE{s} - \QUADE{s_0} & = \QUADF{\tilde s} - \QUADF{s_*}
\end{align*}
and $| \tilde s - s_* | \ge |s-s_0| $.
\end{lemma}

We define for $s^2 \ge \tau_*^2 -
 \tau_{\rm min}^2$
\begin{align*}
 \kappa_s^2 \defn \max \biggl \{ P(f-f^0): \ f \in \calF , \ P(f-f^0)
 + \pen (f) \le \tau_*^2 + s^2 \biggr \} 
 .
\end{align*}

\begin{definition} 
\label{potential.definition} 
We say that the oracle potential holds if
\begin{align}
\label{EqnOraclePotential}
\Gamma \defn \sup_{s >0} \kappa_s /s < \infty .
\end{align}
\end{definition} 

For the shifted version, the counterpart of Condition
\ref{r0.condition} replaces ${\bf E} (\cdot)$ by ${\bf F} (\cdot)$.
\begin{condition} \label{r1.condition}
There is a constant $m_n$, a strictly increasing function $\cal J$
such that the function $\Phi_{\cal J}$ is strictly convex, and such
that the bound
\begin{align}
{\bf F}(s) \le \frac{{\cal J} (s)}{m_n}
\end{align}
holds for all $s \ge 0$.
\end{condition}

When Conditions~\ref{uniformly-bounded.equation}
and~\ref{r1.condition} hold, we define
\begin{align*}
 r_*^2 \defn 2 C^2 \PhiConj_{\cal J} ( 4 K / (m_n C^2)),
\end{align*}
where $\PhiConj_{\cal J}$ is the convex conjugate of $\Phi_{\cal J}$.

\begin{theorem}\label{local2.theorem} Suppose 
that:
\bcar
\item Conditions~\ref{uniformly-bounded.equation}
  and~\ref{r1.condition}, as well as the quadratic curvature condition
  with constant $C$ hold.
\item The shifted mapping $\QUADF{\cdot}$ satisfies the right-sided
  second order margin condition over the interval $( s_* + \underline
  \delta , \sqrt { \tau_{\rm max}^2 - \tau_*^2} ] $ with margin
    function $G$.
\item The oracle potential condition~\eqref{EqnOraclePotential} holds.
\ecar
Then there is a constant $c_0$ depending on $C$, $K$ and $\Gamma$,
such that for all $t>0$ and for a constant $\delta(t) $ such that $G(
\delta (t)) $ is not larger than
\begin{multline*}
 \Gconj \biggl ( c_0 \sqrt {[ t + \log (1+ \sqrt {n(\tau_{\rm max}^2 -
       \tau_*^2 )})]/n} ) \biggr ) \\
+ c_0 \biggl ( (s_*+r_*) \sqrt {[ t + \log (1+ \sqrt {n(\tau_{\rm
        max}^2- \tau_*^2}) )]/n} \\
 + [ t + \log (1+ \sqrt {n(\tau_{\rm max}^2 -\tau_*^2 )} )]/n \biggr )
\end{multline*}
 one has the deviation inequality
\begin{align}
 \PP \biggl ( \sqrt {\tau_{\rm max}^2- \tau_*^2 } \ge \sqrt { \tau^2
   (\hat f) - \tau_*^2} > s_* + \max\{ \delta (t) , \underline
 \delta\} \biggr ) \le \exp[-t] . 
\end{align}
Moreover, if $s_*^2 > \tau_{\rm min}^2 - \tau_*^2 $ and in fact the
two-sided second order margin condition holds for $\QUADF{\cdot}$ in
the range $[\sqrt {\tau_{\rm min}^2 - \tau_*^2} , s_*) \cup ( s_* +
  \underline \delta , \sqrt { \tau_{\rm max}^2 - \tau_*^2} ] $ with
margin function $G$, then one has the concentration inequality
\begin{align*}
\PP \biggl ( |\sqrt { \tau^2(\hat f) - \tau_*^2} - s_* |> \max\{
\delta (t), \underline \delta\} , \ \tau(\hat f) \le \tau_{\rm max}
\biggr ) \le 2 \exp[-t] .
\end{align*}
\end{theorem}


\section{Proofs}
\label{proofs.section}

This section is devoted to the proofs of all our results.
 

\subsection{Proof of Theorem~\ref{direct.theorem}} 
\label{proofs.section.direct}

After some simple algebra, we may write
\begin{align*}
\hat g = \arg \min_{g \in \R^n} \big \{ \frac{1}{2} \tau^2 (g) -
\epsilon^T g / n \big \},
\end{align*}
where $\tau^2 (g) \defn \| g - g^0 \|_n^2 + \pen (g) $. The
function $\tau^2$ is convex, so that it has a sub-differential,
denoted by $\partial \tau^2 (g)$.  With this notation, the minimizing
argument $\hat g$ must satisfy the relation
\begin{equation}  
\label{inclusion.equation}
\epsilon/n \in \partial \tau^2 (\hat g)/2 .
\end{equation}
Due to the strong convexity and coercivity of $g \mapsto \tau^2(g)$,
the inclusion~\eqref{inclusion.equation} always has a unique solution
$\hat g$.
  
We now use a classical fact from convex analysis
(\cite{rockafellar1970convex}): since the function $g \mapsto
\tau^2(g)/2$ is $(1/n)$-strongly convex, the sub-differential mapping
\mbox{$g \mapsto \partial \tau^2 (g)/2$} is $(1/n)$-strongly monotone,
which means that for any pair of vectors $u, u^{\prime} \in \R^n$
\begin{equation}
\label{monotone.equation}
\| u - u^{\prime} \|_n^2 \le (v- v^{\prime} )^T (u- u^{\prime}),
\end{equation}
where $v, v^{\prime} $ denote any members of $\partial \tau^2 (u)/2$
and $\partial \tau^2 (u^{\prime}) /2$ respectively.  

By Borell's theorem~(\citeyear{borell1975brunn}) on the concentration
of Lipschitz functions of Gaussian vectors, it suffices to show that
the mapping $ \epsilon \mapsto \hat m \defn \| \hat g - g^0 \|_n$ is
Lipschitz with parameter $1/\sqrt{n}$.  Let $\epsilon $ and
$\epsilon^{\prime}$ be two realizations of the noise vector, with
corresponding solutions $\hat g$ and $\hat g^{\prime}$, along with
their associated errors $\hat m = \| \hat g - g^0 \|_n$ and $\hat
m^{\prime} = \| \hat g^{\prime} - g^0 \|_n $.  By the triangle
inequality, we have
\begin{align*}
| \hat m - \hat m^{\prime} | = \biggl | \| \hat g - g^0 \|_n - \| \hat
g^{\prime} - g^0 \|_n \biggr | \le \| \hat g - \hat g^{\prime} \|_n,
\end{align*}
so that it suffices to prove that 
\begin{align*}
\| \hat g - \hat g^{\prime} \|_n \le \| \epsilon - \epsilon^{\prime}
\|_n \; = \; \frac{1}{\sqrt{n}} \|\epsilon - \epsilon^\prime \|_2.
\end{align*}

Now consider the pair $u= \hat g$ and $u^{\prime} = \hat g^{\prime} $,
along with the corresponding elements $v= \epsilon/n$ and $v^{\prime}
= \epsilon^{\prime}/n$. Applying the monotone
property~\eqref{monotone.equation} to these pairs yields the
inequality
\begin{align*}
\| \hat g - \hat g^{\prime} \|_n^2 \le ( \epsilon - \epsilon^{\prime}
)^T ( \hat g - \hat g^{\prime} )/n \leq \| \epsilon -
\epsilon^{\prime} \|_n \| \hat g - \hat g^{\prime} \|_n ,
\end{align*}
where the final step follows from the Cauchy-Schwarz inequality.
Cancelling terms completes the proof.
  \hfill $\sqcup \mkern -12mu \sqcap$
%
 

\subsection{Proofs for Section~\ref{theoretical-version.section}}
\label{proof.minimum.lemma}

In this section, we collect the proofs of all results stated in
Section~\ref{theoretical-version.section}.

\paragraph{Proof of Lemma \ref{minimum.lemma}:} 

For any scalar $s$ and $f \in \calF$ such that $\tau(f) \le s$, we
have
\begin{align*}
P_n ( f - f^0) + \pen (f) = \tau^2(f) - (P_n - P) (f^0 - f) \le s^2 -
(P_n - P) (f^0 - f).
\end{align*}
Consequently, we have 
\begin{align*}
\tau^2 (\hat f) - (P_n-P ) (f^0 - \hat f) & = \min_{f \in \calF}
\biggl \{ \tau^2 (f) - (P_n-P ) (f^0 -f) \biggr \} \\
& = \min_{s \ge \tau_{\rm min} } \min_{ \tau(f) \le s } \biggl \{
\tau^2(f) - (P_n - P) (f^0 - f) \biggr \} \\
& \leq \min_{ s \ge \tau_{\rm min} } \min_{ \tau(f) \le s } \biggl \{
s^2 - (P_n - P) (f^0 - f) \biggr \} \\
& = \min_{ s \ge \tau_{\rm min} } \biggl \{ s^2 - \max_{\tau(f) \le s
} (P_n - P) (f^0 - f) \biggr \} \\
& = \min_{s \ge \tau_{\rm min}} \{ s^2 - \hat {\bf E}_n (s) \} .
\end{align*}
On the other hand, for any $f \in \calF$, we have the lower bound
\begin{align*}
\tau^2(f) - (P_n-P) ( f^0 - f) \ge \tau^2 (f) - \max_{\tilde f \in
  \calF_{\tau(f)} } (P_n - P) (f^0 - \tilde f) = \tau^2(f) - \hat {\bf
  E}_n (\tau(f)),
\end{align*}
which implies that
\begin{align*}
 \tau^2(\hat f) - (P_n-P) ( f^0 - \hat f) & \ge \tau^2 (\hat f) - \hat
     {\bf E}_n (\tau(\hat f) ) \; \ge \; \min_{ s \ge \tau_{\rm min}}
     \{ s^2 - \hat {\bf E}_n (s) \} .
\end{align*}
Since the minimizing argument $\hat s = \arg \min_{s \ge \tau_{\rm
    min}} \{ s^2 - \hat {\bf E}_n (s) \}$ is unique by assumption, we
conclude that $\tau (\hat f) = \hat s$, as claimed.  \hfill $\sqcup
\mkern -12mu \sqcap$


\subsection{Proofs for Section~\ref{ERM.section}}
  
We first state and prove an auxiliary lemma that serves as a tool in
the proof of Theorem~\ref{local.theorem}.

\begin{lemma} \label{margin.lemma} 
Let $G$ be a real-valued function with convex conjugate $\Gconj$.
Then all for positive scalars $a$, $b$ and $c$ such that $G(a) \ge
\Gconj(2b) + 2c$, we have \mbox{$G(a) - ab - c \ge 0$.}
\end{lemma}
  
{\bf Proof of~Lemma \ref{margin.lemma}.} By the Fenchel-Young
inequality, we have
\begin{align*}
ab = a \frac{2 b}{2} \le \frac{G(a)}{2} + \frac{\Gconj( 2b)}{2},
\end{align*}
and consequently,
\begin{align*}
G(a) - ab - c \ge G(a) - \frac{G(a)}{2} - \frac{\Gconj(2b)}{2} - c =
\frac{G(a)}{2} - \frac{\Gconj(2b)}{2}- c \ge 0 .
\end{align*}
\hfill $\sqcup \mkern -12mu \sqcap$

We now turn the proof of the main theorem.
  
{\bf Proof of Theorem~\ref{local.theorem}.}
Assume first that the right-sided second order margin condition holds
with $\underline \delta =0$.  Let $t>0$ be arbitrary, and define
\begin{align*}
\bar z(t) \defn 2C s_0 \sqrt {\frac{t}{n}} + r_0 \sqrt {\frac{t}{n}} +
\frac{2 K t}{3 n} , \quad \mbox{and} \quad \underline z(t) \defn 2 C
s_0 \sqrt {\frac{t}{n}} + r_0 \sqrt {\frac{t}{n}} + \frac{Kt}{n} .
\end{align*}

Our strategy is to apply a ``peeling argument'' so as to transition
from the fixed $s$-result of Theorem \ref{Klein-Rio.theorem} to a
result that holds uniformly in $s$.  For a parameter $\epsilon > 0$ to
be chosen later, define the intervals
\begin{align*}
 I_j \defn ((j-1) \epsilon + \delta + s_0, s_0 + \delta + j \epsilon
 ], \qquad \mbox{for $j = 1, \ldots, J$,}
\end{align*}
where $J \defn \lceil \frac{\tau_{\rm max}}{\epsilon} \rceil$, as well
as the associated probabilities
\begin{align*}
 \PP_j = \PP \biggl ( \exists \ s \in I_j \mbox{ such that }  s^2 -
 \hat {\bf E}_n (s) \le s_0^2 - {\bf E} (s_0) + \bar z(t) \biggr )
\end{align*}
Then for a parameter $\delta >0$ to be chosen later (and leading to
the $\delta(t) $ in the theorem statement), we have
\begin{multline*}
\PP \Big ( s_0 + \delta < \hat s \le \tau_{\rm max} , \ \hat {\bf E}_n
(s_0) > {\bf E}(s_0) + \bar z(t) \Big ) \\ \leq \PP \biggl ( \exists s
\in ( s_0 + \delta, \tau_{\rm max}] \, \mbox{ such that } \, \ s^2 -
  \hat {\bf E}_n (s) \le s_0^2 - {\bf E} (s_0) + \bar z(t) \biggr ) \;
  \leq \sum_{j=1}^{J} \PP_j.
\end{multline*}
For each index $j$ and for all $s \in I_j$, we have
\begin{align*}
s^2 - \hat {\bf E}_n (s) \ge ( (j-1) \epsilon + \delta + s_0 )^2 -
\hat {\bf E}_n ( s_0 + \delta + j \epsilon ).
\end{align*}
Moreover, for all $u>0$, we have  by Theorem \ref{Klein-Rio.theorem}
\begin{align*}
\hat {\bf E}_n (s_0 + \delta + j \epsilon) \le {\bf E}(s_0 + \delta +
j \epsilon ) - 2C (\delta + j \epsilon ) \sqrt {u/n}- \underline z(u)
\end{align*}
with probability at least $1- \exp[-u]$. Furthermore, by the one-sided
form of second order margin condition, we have the lower bound
\begin{multline*} 
(s_0 + \delta + (j -1) \epsilon )^2 - {\bf E}(s_0 + \delta + j
  \epsilon ) \ge s_0^2 - {\bf E}(s_0) + G( \delta + j \epsilon) \\
+ (s_0 + \delta (j -1) \epsilon )^2 - (s_0 + \delta +j \epsilon )^2
\\ 
= s_0^2 - {\bf E}(s_0) + G( \delta + j \epsilon) - 2\epsilon ( s_0 +
\delta + j \epsilon ) + \epsilon^2 .
\end{multline*}
Putting together the pieces, for all $s \in I_j$, we have
\begin{multline*}
 s^2 - \hat {\bf E}_n (s) \ge s_0^2 - {\bf E}(s_0) + G( \delta + j
 \epsilon) - 2 \Big(C \sqrt {\frac{u}{n}} +\epsilon \Big) \: ( \delta
 + j \epsilon ) - 2\epsilon s_0 + \epsilon^2 - \underline z(u),
\end{multline*}
with probability at least $1- \exp[-u]$.  We now apply
Lemma~\ref{margin.lemma} with the choices \mbox{$a \defn \delta + j
  \epsilon$,} \mbox{$b \defn 2(C \sqrt {u/n} +\epsilon ) $,} and
\mbox{$c\defn 2\epsilon s_0 - \epsilon^2 + \underline z(u) + \bar z(t)
  $.}  In order to be able to do so, we require that 
\begin{align}
\label{EqnRequire}
G(\delta) \ge \Gconj \biggl ( 4(C \sqrt {u/n} +\epsilon ) \biggr ) + 2
\biggl ( 2\epsilon s_0 - \epsilon^2 + \underline z(u) + \bar z(t)
\biggr ) .
\end{align}

We now settle the choice of $\epsilon$ and $u$. Taking $\epsilon =
1/\sqrt n $, we are then guaranteed that
\begin{align*}
J \leq 1+ { \tau_{\rm max} \over \epsilon } = 1+ \sqrt {n\tau_{\rm
    max}^2}.
\end{align*}
Moreover, recalling the arbitrary $t > 0$ introduced at the start of
the proof, we set $u= t + \log (1+ \sqrt {n\tau_{\rm max}^2} ) $, and
then the condition~\eqref{EqnRequire} on $\delta:= \delta(t)$ becomes
\begin{multline*}
G(\delta (t)) \defn \Gconj \biggl ( 4(C \sqrt {[ t + \log (1+ \sqrt
    {n\tau_{\rm max}^2} )]/n} +1/ \sqrt n ) \biggr ) \\
+ 2 \biggl ( 2s_0/\sqrt n -1/n + \underline z([ t + \log (1+ \sqrt
  {n\tau_{\rm max}^2} )]) + \bar z(t) \biggr ). 
\end{multline*}
We are then guaranteed that for each $j \in \{1 ,
\ldots , J\}$, with probability at least $1- \exp[-(t+ \log (1+ \sqrt
  {n \tau_{\rm max}^2} ) ] $, for all $s \in I_j$ it holds that
\begin{align*}
 s^2 - \hat {\bf E}_n (s) \ge s_0^2 - {\bf E}(s_0) + \bar z(t) . 
\end{align*}
It follows that $\PP_j \le \exp \Big \{ -[t+ \log (1+ \sqrt
  {n\tau_{\rm max}^2})] \Big \}$, and hence
\begin{align*}
\sum_{j=1}^J \PP_j \le J \exp\biggl [-[t+ \log (1+ \sqrt {n\tau_{\rm
        max}^2} ] \biggr ] \le \exp [-t] . 
\end{align*}
One easily verifies that for some constant $c_0$ depending on $C$ and
$K$, we have
\begin{multline*}
G(\delta (t)) \le \Gconj \biggl ( c_0 \sqrt {[ t + \log (1+ \sqrt
    {n\tau_{\rm max}^2} )]/n} ) \biggr ) \\+ c_0 \biggl ( (s_0+r_0)
\sqrt {[ t + \log (1+ \sqrt {n\tau_{\rm max}^2} )]/n} + [ t + \log (1+
  \sqrt {n\tau_{\rm max}^2} )]/n \biggr ).
\end{multline*}
Overall, we have established that $\PP \big ( \hat s > s_0 + \delta
(t) \big ) \le 2 \exp[-t]$.

In stating the bound in the theorem, we removed the pre-factor of 2
for cosmetic reasons. This can be done by replacing $t $ by $t+ \log
2$.  In order to prove the lower bound, one may follow the same
argument, instead using the left-sided version of the second order
margin condition.

In our argument thus far, we assumed $\underline \delta =0$.  If the
second order margin condition only holds at distance $\underline
\delta >0$, it is clear that one simply can take the maximum of
$\delta (t)$ and $\underline \delta$ in the bounds.  \hfill $\sqcup
\mkern -12mu \sqcap$


\subsection{Proofs for Section~\ref{quadratic.section}}

{\bf Proof of Lemma \ref {concaveq.lemma}.} Let us introduce the
shorthand notation $\tilde s \defn s^{2/q}$ and $\tilde s_0 =
s_0^{2/q}$, let $\tilde s_1\ge \tau_{\rm min}^{2/q} $ and $\tilde s_2
\ge \tau_{\rm min}^{2/q}$ be arbitrary, and define
\begin{align*}
\hat f_1 & \defn f_{\hat g_1} \defn \arg \max_{ \tau^{2/q} (f) \le
  \tilde s_1 } (P_n - P) (f^0- f), \quad \mbox{as well as} \\
\hat f_2 & \defn f_{\hat g_2} \defn \arg \max_{ \tau^{2/q} (f) \le
  \tilde s_2 } (P_n - P) (f^0- f ).
\end{align*}
With these choices, we have
\begin{align*}
\tau^{2/q} (f_{t \hat g_1+ (1- t) \hat g_2} ) \le t \tau^{2/q} (\hat
f_1) + (1-t) \tau^{2/q} (\hat f_2) \le t \tilde s_1 + (1- t) \tilde
s_2 \qquad \mbox{for all $t \in [0,1]$.}
\end{align*}
Moreover, we have the lower bound
\begin{align*}
 \hat {\bf E}_{n } (( t\tilde s_1 + (1- t) \tilde s_2)^{q/2})\ge (P_n
 - P) ( f^0 - f_{t \hat g_1 + (1-t) \hat g_2} ) \\
= t(P_n-P) (f^0 - \hat f_1) + (1-t) (P_n-P) (f^0 - \hat f_2) \\
= t \hat {\bf E}_{n} (\tilde s_1^{q/2}) + (1-t) \hat {\bf E}_{n }
(\tilde s_2^{q/2} ) .
\end{align*}
Taking expectations yields the lower bound
\begin{align*}
 {\bf E} (( t\tilde s_1 + (1- t) \tilde s_2)^{q/2}) & \ge t {\bf E}
 (\tilde s_1^{q/2}) + (1-t) {\bf E} (\tilde s_2^{q/2} ) .
\end{align*}
Using the fact that
\begin{align*}
[( t\tilde s_1 + (1- t) \tilde s_2)^q - {\bf E} ( ( t\tilde s_1 + (1-
  t) \tilde s_2)^{q/2}) ]- [\tilde s_0^q - {\bf E} ( \tilde
  s_0^{q/2})] \ge 0,
\end{align*}
we have
\begin{align*}
 [ \tilde s^q - {\bf E} ( \tilde s^{q/2} ) ]- [ \tilde s_0^q - {\bf E}
   ( \tilde s_0^{q/2} ) ] \ge \tilde s^q - \tilde s_0^q + {1 \over t}
 { \tilde s_0^q - ( t \tilde s + (1- t) \tilde s_0)^q } .
\end{align*}
Taking $t \downarrow 0$ then gives
\begin{align*}
[ \tilde s^q - {\bf E}( \tilde s^{q/2} ) ]- [ \tilde s_0^q - {\bf E} (
  \tilde s_0^{q/2} ) ] & \ge \tilde s^q - \tilde s_0^q - q\tilde
s_0^{q-1} (\tilde s - \tilde s_0) \\
& \ge q (q-1)(M+1)^{-{2(2-q) \over q} } \tilde s_0^{-(2-q)} (\tilde s
- \tilde s_0)^2 / 2,
\end{align*}
valid when $q \in (1,2]$ and $\tilde s > (M+1)^{2/q} \tilde s_0$.
  Furthermore, for $1 < q \le 2$ we get for some $s_0 < \bar s \le s $
\begin{align*}
 \tilde s - \tilde s_0 = s^{2/q} - s_0 ^{2/q} = 2 \bar s^{{2-q \over
     q} } (s- s_0)/q \ge {2 } s_0^{{2-q \over q} } (s- s_0) /q. 
\end{align*}
Consequently, for all $q \in (1,2]$, we have
\begin{align*} 
[ \tilde s^q - {\bf E} ( \tilde s^{q/2} ) ]- [ \tilde s_0^q - {\bf E}
  ( \tilde s_0^{q/2} ) ] \ge 2 (q-1) (M+1)^{-{2(2-q) \over q} }
(s-s_0)^2 /q,
\end{align*}
as claimed.  \hfill $\sqcup \mkern -12mu \sqcap$
  

{\bf Proof of Lemma \ref{approx-concave.lemma}.}  For all $ s \ge
\tau_{\rm min}$, we have
\begin{align*}
 {\bf E}^{\tau} (As ) \ge {\bf E}^{\varsigma} (A^{1/2} s) \ge {\bf
   E}^{\tau} (s) .
\end{align*}
For $t s + (1- t) s_0\ge \tau_{\rm min}$, the concavity of ${\bf
  E}^{\varsigma}$ yields
\begin{align*}
{\bf E}^{\tau} \big ( A ( ts + (1- t) s_0) \big ) & \ge {\bf
  E}^{\varsigma} \big ( A^{1/2}( ts + (1- t) s_0) \big ) \\
 & \ge t {\bf E}^{\varsigma} ( A^{1/2}s) + (1- t) {\bf E}^{\varsigma}
(A^{1/2}s_0) \\
& \ge t {\bf E}^{\tau} (s) + (1- t) {\bf E}^{\tau} (s_0) \\
& = t({\bf E}(s) - {\bf E}(s_0) ) + {\bf E}(s_0).
\end{align*}
Since $s_0$ is the minimizer of $s^2 - {\bf E}(s)$, we have
\begin{align*}
A^2 \big ( ts + (1-t) s_0 \big )^2 - {\bf E} \big ( A^2 ( ts + (1-t)
s_0) \big ) \ge s_0^2 - {\bf E}(s_0),
\end{align*}
and consequently
\begin{align*}
{\bf E} \big ( A ( ts + (1-t) s_0) \big ) \le A^2 \big ( ts + (1-t)
s_0 \big )^2 -s_0^2 + {\bf E}(s_0).
\end{align*}
It follows that
\begin{align*}
t ({\bf E}(s) - {\bf E}(s_0) ) & \le A^2 \big ( ts + (1-t) s_0 \big
)^2 -s_0^2 \\
& = A^2 ( t^2 s^2 + s_0^2 -2ts_0^2 + t^2 s_0^2 + 2t s s_0 - 2 t^2 s
s_0) - s_0^2 \\
& = A^2 ( t^2 (s-s_0)^2 + s_0^2 -2t s_0^2 + 2t s s_0 ) - s_0^2 \\
 & = t^2 (s-s_0)^2 -2t s_0^2 + 2t s s_0 + \epsilon \biggr ( t^2 (s-
s_0)^2 -2t s_0^2+ 2t s s_0 \biggr ) + \epsilon s_0^2 .
 \end{align*}
Putting together the pieces, we have
\begin{multline*}
 [s^2 - {\bf E}(s)] - [s_0^2 - {\bf E}(s_0)] = s^2 - s_0^2 - { t (
   {\bf E}(s) - {\bf E}(s_0)) \over t } \\
\ge s^2 - s_0^2 - t (s-s_0)^2
 +2 s_0^2 - 2s s_0 - \epsilon \biggr ( t (s- s_0)^2 -2 s_0^2+ 2 s s_0
 \biggr ) - \epsilon s_0^2/t \\
 = (s-s_0)^2 - t (s-s_0)^2 -\epsilon \biggr ( t (s- s_0)^2 +2
 s_0(s-s_0) \biggr ) - \epsilon s_0^2 /t .
\end{multline*}

We now choose $t = \sqrt \epsilon$ to find that
\begin{multline*}
[s^2 - {\bf E}(s)] - [s_0^2 - {\bf E}(s_0)] \ge (s-s_0)^2 - \sqrt
\epsilon (s-s_0)^2 - \epsilon \biggr ( \sqrt \epsilon (s- s_0)^2 +2
s_0(s-s_0) \biggr ) - \sqrt {\epsilon} s_0^2 .
\end{multline*}
Next we use our assumption that $s \le \tau_{\rm max}= (M+1)s_0$. We
then get
\begin{align*}
[s^2 - {\bf E}(s)] - [s_0^2 - {\bf E}(s_0)] & \ge \biggl (1- \sqrt
\epsilon (1+ \epsilon)\biggr ) (s-s_0)^2 - \epsilon \biggr ( 2 s_0^2M
\biggr ) - \sqrt {\epsilon} s_0^2 \\
& = (1- \sqrt \epsilon (1+ \epsilon) ) (s-s_0)^2- \sqrt \epsilon ( 2
\sqrt \epsilon M + 1) s_0^2 \\
& \ge (s-s_0)^2 / 2 - \sqrt \epsilon ( 2 \sqrt \epsilon M + 1) s_0^2 .
\end{align*}
Now we take $| s- s_0| \ge 2 \big [ \sqrt \epsilon ( 2 \sqrt \epsilon
  M + 1) \big ]^{1/2} s_0$, and conclude that
\begin{align*}
[s^2 - {\bf E}(s)] - [s_0^2 - {\bf E}(s_0)] \ge (s-s_0)^2 / 4 .
\end{align*}
\hfill $\sqcup \mkern -12mu \sqcap$


{\bf Proof of Lemma~\ref{concave.lemma}.}  Let $s_1\ge \varsigma_{\rm
  min}$ and $s_2 \ge \varsigma_{\rm min}$ be arbitrary, and define
\begin{align*}
\hat f_1 \defn f_{\hat g_1} \defn \arg \max_{ \varsigma (f) \le s_1 }
(P_n - P) (f^0- f), \quad \mbox{and} \quad \hat f_2 \defn f_{\hat g_2}
\defn \arg \max_{ \varsigma (f) \le s_2 } (P_n - P) (f^0- f ).
\end{align*}
For all $t \in [0,1]$, we have
\begin{align*}
\varsigma (f_{t \hat g_1+ (1- t) \hat g_2} ) \le t \varsigma (\hat
f_1) + (1-t) \varsigma (\hat f_2) \le t s_1 + (1- t) s_2.
\end{align*}
In addition, we have
\begin{align*}
 \hat {\bf E}_{n}^{ \varsigma} ( ts_1 + (1- t) s_2) & \ge (P_n - P) (
 f^0 - f_{t \hat g_1 + (1-t) \hat g_2} ) \\
& = t(P_n-P) (f^0 - \hat f_1) + (1-t) (P_n-P) (f^0 - \hat f_2)\\
& = t \hat {\bf E}_{n }^{ \varsigma} (s_1) + (1-t) \hat {\bf E}_{n }^{
   \varsigma} (s_2),
\end{align*}
which completes the proof.  \hfill $\sqcup \mkern -12mu \sqcap$


\subsection{Proofs for Section~\ref{exponential-families.section}}
  
{\bf Proof of Lemma~\ref{exponential-family.lemma}.} Throughout this
proof, we let $0 \le \tilde t \le t$ be some intermediate point, not
the same at each appearance.  The function $h(t) = \dfun(t g)$ is
infinitely differentiable with $h(0) = 0$ and $h'(0) = 0$, so a
second-order Taylor series expansion yields $\dfun (t g) =  {1
  \over 2} t^2 h''(t)$.  Consequently, it suffices to show that
\begin{align}
\label{EqnSecondFinal}
h''(t) = Pg^2 ( 1+ {\mathcal O} (t)).
\end{align}
Computing derivatives, we have
\begin{align*}
h'(t) & = \big [ P \exp [t g] \big ]^{-1} P (\exp [tg ] g), \quad
\mbox{and} \\
h''(t) & = \big [ P \exp [t g ] \big ]^{-1} P (\exp [tg] g^2) - \big
\{ \big [ P \exp [t g] \big ]^{-1} P (\exp [tg ] g) \big \} ^2 .
\end{align*}

By a Taylor series expansion of the exponential, we have
\begin{align*} 
P \exp [tg] & = 1 + t Pg + t^2 P (\exp [ \tilde t g ] g^2) /2 = 1+
{\mathcal O} (t^2 ) Pg^2= 1 + {\mathcal O} (t^2), \quad \mbox{and} \\
 P \exp [tg] g & = t P g^2 + {t^2 \over 2}  P (\exp [ \tilde t g ] g^3) = t Pg^2
 + {\mathcal O} (t^2) Pg^2 = t Pg^2 ( 1+ {\mathcal O} (t)).
\end{align*}
Combining the pieces, we find that
\begin{align*}
\big [ P \exp [t g] \big ]^{-1} P (\exp [tg ] g ) = [ 1+ {\mathcal O}
  (t^2 ) ]^{-1} \big [ tPg^2 (1+ {\mathcal O}(t)) \big ] = t Pg^2 (1+
     {\mathcal O} (t) ) .
\end{align*}
It follows that
\begin{align*}
\big \{ \big [ P \exp [t g] \big ]^{-1} P (\exp [tg ] g )\big \}^2 =
t^2 (Pg^2 )^2 (1+ {\mathcal O} (t) ) = {\mathcal O} (t^2 Pg^2 ).
\end{align*}
But $P (\exp [tg] g^2) = Pg^2 + t P( \exp[\tilde t g ] g^3 = Pg^2 (1+
{\mathcal O} (t) ) $, and hence the bound~\eqref{EqnSecondFinal}
follows, which completes the proof.

\hfill $\sqcup \mkern -12mu
 \sqcap$


\subsection{Proofs for Section~\ref{concentration.section}}
 
{\bf Proof of Lemma \ref{kleinrio-truncated.lemma}.}  

For each
$t>0$, we have
\begin{align*}
(P_n - P) (f^0- f) \le (P_n - P) (f^0-f) {\rm l} \{ F \le t \} + (P_n
  - P) F {\rm l} \{ F > t \} + 2P F {\rm l} \{ F > t \},
\end{align*}
and also
\begin{align*}
(P_n - P) (f^0- f) \ge (P_n - P) (f^0-f) {\rm l} \{ F \le t \} - (P_n
  - P) F {\rm l} \{ F > t \} - 2P F {\rm l} \{ F > t \} .
\end{align*}
Taking $t$ here equal to $t _0: = C_F \sqrt {\log n }$ (and assuming
$t_0>1$) we see that
$$ 2 P F {\rm l }\{ F > t_0 \} \le 2 P F^2 {\rm l} \{ F >t_0 \}  \le 2 c_F^2 /n .$$
Moreover, for all $t>0$,  with probability at least $1-1/t^2$
$$ | (P_n - P) F {\rm l} \{ F > t_0 \} | \le t (P (F^2 {\rm l } \{ F
>t_0 \} )^{1/2} / \sqrt n \le t c_F / n . $$ Write the truncated
versions as
\begin{align*}
\hat {\bf E}_n^{\rm trunc} (s) & \defn \max_{f \in \calF_s } (P_n - P)
(f^0-f) {\rm l} \{ F \le t \} , \quad \mbox{and} \\
{\bf E}^{\rm trunc} (s) & \defn \EE \biggl ( \max_{f \in \calF_s }
(P_n - P) (f^0-f) {\rm l} \{ F \le t \} \biggr ) .
\end{align*}
Then $\big | {\bf E} (s) - {\bf E}^{\rm trunc} (s) \big | \le 2c_F^2
/n$, and moreover, with probability at least \mbox{$1- 1/t^2$,} we
have
\begin{align*}
\Big | \hat {\bf E}_n (s) - \hat {\bf E}_n^{\rm trunc} (s) \Big | &
\le c_F 2 (c_F + t ) /n .
\end{align*}

Now Theorem~\ref{Klein-Rio.theorem} ensures that, for all $t \geq 0$,
\begin{align*}
\hat {\bf E}_n^{\rm trunc} (s) & \ge {\bf E}^{\rm trunc}(s) - \sqrt {
  8 C_F \sqrt {\log n} {\bf E}^{\rm trunc} (s) + 2 \sigma_s^2 } \sqrt
     {t/n} - C_F t \sqrt {\log n} / n , \quad \mbox{and} \\
\hat {\bf E}_n^{\rm trunc} (s) & \le {\bf E}^{\rm trunc} (s) + \sqrt {
  8 C_F \sqrt {\log n} {\bf E}^{\rm trunc}(s) + 2 \sigma_s^2 } \sqrt
     {t/n} + 2C_F t \sqrt {\log n} / (3n) ,
\end{align*}
where each bound holds with probability at least $1- \exp[ -t]$,
\hfill $\sqcup \mkern -12mu \sqcap$

\vspace*{.05in}

{\bf Proof of Lemma \ref{Klein-Rio.lemma}.}  By Theorem
\ref{Klein-Rio.theorem} of Klein and Rio, for all $t \ge 0$, with
probability at least $1- \exp[ -t]$,
$$ \hat {\bf E}_n (s) \ge {\bf E}(s) - \sqrt { 8 K {\bf E}(s) + 2 \sigma_s^2 } \sqrt {t/n} - K t/ n  $$
But 
\begin{align*}
8K {\bf E}(s) \le 8K {\cal J} (s) / m_n \le 2 C^2 s^2 + 2 C^2 \PhiConj
( 4 K / ( m_n C^2 ))= 2 C^2 s^2 + r_0^2 .
\end{align*}
Moreover, we have $ \sigma_s^2 \le C^2 s^2 $ by the quadratic curvature
condition, and hence
\begin{align*}
 \sqrt { 8 K {\bf E}(s) + 2 \sigma_s^2 } \sqrt {t/n} \le \sqrt { 4 C^2
   s^2 + r_0^2 } \sqrt {t/n} \le 2 C s \sqrt {t/n} + r_0 \sqrt {t/n}.
\end{align*}
\hfill $\sqcup \mkern -12mu \sqcap$

 \subsection{Proofs for Section \ref{shifted.section}}
 
 {\bf Proof of Lemma \ref{shifted.lemma}.} 
Since ${\bf F}(\tilde s) = {\bf E} ( s  )$,
it follows that
\begin{align*}
\tilde s^2 - {\bf F}(\tilde s) = s^2 - {\bf E}(s) - \tau_*^2 . 
\end{align*}
We also have
\begin{align*}
|s-s_0| = \biggl | \sqrt {\tilde s^2 + \tau_*^2 } - \sqrt { s_*^2 +
  \tau_*^2 } \biggr | & = | \tilde s - s_* | { \tilde s + s_* \over
  \sqrt { \tilde s^2 + \tau_*^2} + \sqrt { \tilde s_*^2 + \tau_*^2} }
\\
& \le | \tilde s - s_* |
\end{align*}
 \hfill $\sqcup \mkern -12mu \sqcap$


{\bf Proof of Theorem \ref{local2.theorem}.} By the oracle potential
(see Definition~\ref{potential.definition}), we know that for all
scalars $s$ such that \mbox{$s^2 \ge \tau_*^2 - \tau_{\rm min}^2$,} we
have 
\begin{align*}
\sup_{f \in \calF : \tau^2 ( f) \le \tau_*^2 + s^2 } P(f-f_0)
\le \Gamma^2 s^2,
\end{align*}
and hence $\sup \limits_{f \in \calF , : \, \tau^2 ( f) \le \tau_{*}^2
  + s^2 } \sigma^2 (f-f^0) \le \Gamma^2 C^2 s^2$.  We can thus proceed
along the same lines as in the proof of Theorem~\ref{local.theorem},
replacing $C$ by $\Gamma C$.

\hfill $\sqcup \mkern -12mu \sqcap $

\bibliographystyle{plainnat}
\bibliography{reference}

\begin{thebibliography}{13}
\providecommand{\natexlab}[1]{#1}
\providecommand{\url}[1]{\texttt{#1}}
\expandafter\ifx\csname urlstyle\endcsname\relax
  \providecommand{\doi}[1]{doi: #1}\else
  \providecommand{\doi}{doi: \begingroup \urlstyle{rm}\Url}\fi

\bibitem[Borell(1975)]{borell1975brunn}
C.~Borell.
\newblock The {B}runn-{M}inkowski inequality in {G}auss space.
\newblock \emph{Inventiones Mathematicae}, 30\penalty0 (2):\penalty0 207--216,
  1975.

\bibitem[Boucheron and Massart(2011)]{boucheron2011high}
S.~Boucheron and P.~Massart.
\newblock A high-dimensional {W}ilks phenomenon.
\newblock \emph{Probability Theory and Related Fields}, 150\penalty0
  (3-4):\penalty0 405--433, 2011.

\bibitem[Boucheron et~al.(2013)Boucheron, Lugosi, and
  Massart]{boucheron2013concentration}
S.~Boucheron, G.~Lugosi, and P.~Massart.
\newblock \emph{Concentration inequalities: A nonasymptotic theory of
  independence}.
\newblock OUP Oxford, 2013.

\bibitem[Chatterjee(2014)]{chatterjee2014new}
S.~Chatterjee.
\newblock A new perspective on least squares under convex constraint.
\newblock \emph{The Annals of Statistics}, 42\penalty0 (6):\penalty0
  2340--2381, 2014.

\bibitem[Klein(2002)]{klein2002inegalite}
T.~Klein.
\newblock Une in{\'e}galit{\'e} de concentration {\`a} gauche pour les
  processus empiriques.
\newblock \emph{Comptes Rendus Mathematique}, 334\penalty0 (6):\penalty0
  501--504, 2002.

\bibitem[Klein and Rio(2005)]{klein2005concentration}
T.~Klein and E.~Rio.
\newblock Concentration around the mean for maxima of empirical processes.
\newblock \emph{The Annals of Probability}, 33\penalty0 (3):\penalty0
  1060--1077, 2005.

\bibitem[Koltchinskii(2011)]{koltchinskii2011oracle}
V.~Koltchinskii.
\newblock \emph{Oracle Inequalities in Empirical Risk Minimization and Sparse
  Recovery Problems: Ecole dÕEt{\'e} de Probabilit{\'e}s de Saint-Flour
  XXXVIII-2008}, volume~38.
\newblock Springer Science \& Business Media, 2011.

\bibitem[Ledoux(2001)]{ledoux2001concentration}
M.~Ledoux.
\newblock \emph{The concentration of measure phenomenon}, volume~89.
\newblock American Mathematical Society, 2001.

\bibitem[Massart(2000)]{massart2000some}
P.~Massart.
\newblock Some applications of concentration inequalities to statistics.
\newblock In \emph{Annales de la Facult{\'e} des sciences de Toulouse:
  Math{\'e}matiques}, volume~9, pages 245--303, 2000.

\bibitem[Muro and van~de Geer(2015)]{muro15}
A.~Muro and S.~van~de Geer.
\newblock Concentration behavior of the penalized least squares estimator,
  2015.
\newblock arXiv:1511.08698.

\bibitem[Rockafellar(1970)]{rockafellar1970convex}
R.T. Rockafellar.
\newblock \emph{Convex {A}nalysis}.
\newblock Princeton University Press, 1970.

\bibitem[Saumard(2012)]{saumard2012optimal}
A.~Saumard.
\newblock Optimal upper and lower bounds for the true and empirical excess
  risks in heteroscedastic least-squares regression.
\newblock \emph{Electronic Journal of Statistics}, 6:\penalty0 579--655, 2012.

\bibitem[Talagrand(1995)]{Talagrand:95}
M.~Talagrand.
\newblock {Concentration of measure and isoperimetric inequalities in product
  spaces}.
\newblock \emph{Publications Math\'ematiques de l'IHES}, 81:\penalty0 73--205,
  1995.

\end{thebibliography}
\end{document}